\documentclass[11pt]{amsart}
\usepackage{amsmath,amssymb,latexsym,esint,cite,mathrsfs}
\usepackage{verbatim,wasysym}
\usepackage[left=2.6cm,right=2.6cm,top=2.9cm,bottom=2.9cm]{geometry}
\usepackage{xcolor}

\usepackage{graphicx}
\usepackage{subfigure}

\usepackage{graphicx,epic,eepic}

\allowdisplaybreaks

\begin{document}

\title[Normalized solutions for Schr\"{o}dinger systems]
{Normalized solutions for Schr\"{o}dinger systems in dimension two}

\author[S. Deng]{Shengbing Deng}
\address{\noindent S. Deng-School of Mathematics and Statistics, Southwest University,
Chongqing 400715, People's Republic of China}\email{shbdeng@swu.edu.cn}

\author[J. Yu]{Junwei Yu}
\address{\noindent J. Yu-School of Mathematics and Statistics, Southwest University,
Chongqing 400715, People's Republic of China.}\email{JwYumaths@163.com}

\maketitle

\maketitle
\numberwithin{equation}{section}
\newtheorem{theorem}{Theorem}[section]
\newtheorem{lemma}[theorem]{Lemma}
\newtheorem{corollary}[theorem]{Corollary}
\newtheorem{definition}[theorem]{Definition}
\newtheorem{proposition}[theorem]{Proposition}
\newtheorem{remark}[theorem]{Remark}
\allowdisplaybreaks

\maketitle

\noindent {\bf Abstract}: In this paper, we study the existence of normalized solutions to the following nonlinear Schr\"{o}dinger systems with exponential growth
\begin{align*}
 \left\{
\begin{aligned}
&-\Delta u+\lambda_{1}u=H_{u}(u,v), \quad
\quad
\hbox{in }\mathbb{R}^{2},\\
&-\Delta v+\lambda_{2} v=H_{v}(u,v), \quad
\quad
\hbox{in }\mathbb{R}^{2},\\
&\int_{\mathbb{R}^{2}}|u|^{2}dx=a^{2},\quad \int_{\mathbb{R}^{2}}|v|^{2}dx=b^{2},
\end{aligned}
\right.
\end{align*}
where $a,b>0$ are prescribed, $\lambda_{1},\lambda_{2}\in \mathbb{R}$ and the functions $H_{u},H_{v}$ are partial derivatives of a Carath\'{e}odory function $H$ with $H_{u},H_{v}$ have exponential growth in $\mathbb{R}^{2}$. Our main results are totally new for Schr\"{o}dinger systems in $\mathbb{R}^{2}$. Using the Pohozaev manifold and variational methods, we establish the existence of normalized solutions to the above problem.

\vspace{3mm} \noindent {\bf Keywords}: Normalized solution; Nonlinear Schr\"odinger systems; Variational methods; Critical exponential growth; Trudinger-Moser inequality.

\vspace{3mm} \noindent {\bf 2020  Mathematics Subject Classification.} Primary 35J47;  Secondly 35B33, 35J50.

\vspace{3mm}

\maketitle

\section{{\bfseries Introduction}}\label{introduction}

This paper concerns the existence of normalized solutions to the following nonlinear Schr\"{o}dinger systems with exponential growth
\begin{equation}\label{aa}
 \left\{
\begin{aligned}
&-\Delta u+\lambda_{1}u=H_{u}(u,v), \quad
\quad
\hbox{in }\mathbb{R}^{2},\\
&-\Delta v+\lambda_{2} v=H_{v}(u,v), \quad
\quad
\hbox{in }\mathbb{R}^{2},\\
&\int_{\mathbb{R}^{2}}|u|^{2}dx=a^{2},\quad \int_{\mathbb{R}^{2}}|v|^{2}dx=b^{2},
\end{aligned}
\right.
\end{equation}
where $a,b>0$ are prescribed, $\lambda_{1},\lambda_{2}\in \mathbb{R}$ and the functions $H_{u},H_{v}$ are partial derivatives of a Carath\'{e}odory function $H$ with $H_{u},H_{v}$ have exponential growth in $\mathbb{R}^{2}$.

We say that a function $h(s)$ has subcritical exponential growth if for all $\gamma>0$ we have
\begin{equation}\label{ab}
    \begin{aligned}\displaystyle
    \lim\limits_{|s|\rightarrow +\infty}\frac{|h(s)|}{e^{\gamma s^{2}}}=0
    \end{aligned}
\end{equation}
and $h(s)$ has critical exponential growth if there exists $\gamma_{0}>0$ such that
\begin{equation}\label{ac}
    \begin{aligned}\displaystyle
    \lim\limits_{|s|\rightarrow +\infty}\frac{|h(s)|}{e^{\gamma s^{2}}}=0,\ \ \forall \ \gamma>\gamma_{0},\ \ and \ \ \lim\limits_{|s|\rightarrow +\infty}\frac{|h(s)|}{e^{\gamma s^{2}}}=+ \infty,\ \ \forall \ \gamma<\gamma_{0}.
    \end{aligned}
\end{equation}
This notion of criticality is motivated by the inequality of Trudinger and Moser; see \cite{Moser,Trudinger}.

The problem under consideration arises from the time-dependent of coupled nonlinear Schr\"{o}dinger systems
\begin{equation}\label{abb}
 \left\{
\begin{aligned}
&i \frac{\partial}{\partial t}\Phi_{1}=\Delta \Phi_{1}-V(x)\Phi_{1}+H_{1}(\Phi_{1},\Phi_{2}), \quad
\quad
\hbox{in }\mathbb{R}^{N}\times \mathbb{R},\\
&i \frac{\partial}{\partial t}\Phi_{2}=\Delta \Phi_{2}-V(x)\Phi_{2}+H_{2}(\Phi_{1},\Phi_{2}), \quad
\quad
\hbox{in }\mathbb{R}^{N}\times \mathbb{R},
\end{aligned}
\right.
\end{equation}
where $H_{1},H_{2}$ are partial derivatives of a Carath\'{e}odory function $H$. An important, and of course well known, feature of \eqref{abb} is conservation of mass with the $L^{2}$-norms $|\Phi_{1}(\cdot,t)|_{2},$ $|\Phi_{2}(\cdot,t)|_{2}$ of solutions are independent of $t\in \mathbb{R}$.

Suppose $H_{u}(e^{i\lambda_{1} t}u(x),e^{i\lambda_{2} t}v(x))=e^{i\lambda_{1} t} H_{1}(u,v)$ and $H_{v}(e^{i\lambda_{1} t}u(x),e^{i\lambda_{2} t}v(x))=e^{i\lambda_{2} t} H_{2}(u,v)$ hold. The ansatz $\Phi_{1}(x,t)=e^{i\lambda_{1} t}u(x)$ and $\Phi_{2}(x,t)=e^{i\lambda_{2} t}v(x)$ for solitary wave solutions leads to elliptic system
\begin{equation}\label{b}
 \left\{
\begin{aligned}
&-\Delta u+(\lambda_{1}+V(x))u=H_{u}(u,v), \quad
\quad
\hbox{in }\mathbb{R}^{N},\\
&-\Delta v+(\lambda_{2}+V(x)) v=H_{v}(u,v), \quad
\quad
\hbox{in }\mathbb{R}^{N},
\end{aligned}
\right.
\end{equation}
In literature, few people have devoted themselves to the study of coupled systems containing exponential nonlinearities. If $N=2$, $\lambda_{1},\lambda_{2}=0$ are prescribed and $V(x)$ satisfies some suitable conditions, system \eqref{b} has been investigated by some authors with exponential critical nonlinearities; see e.g. \cite{Alba,Albb,Cas,DOJM,DOJMa} and references therein. If $N=1$, do \'{O} et al. \cite{DOGI} studied the nonautonomous fractional Hamiltonian system with critical exponential growth in $\mathbb{R}$ via Linking Theorem and Galerkin approximation procedure. For more results on $(p, q)$-fractional Laplace system, the reader may refer to \cite{Chen,Thin} and references therein.

When the $L^{2}$-norms $|u|_{2}$, $|v|_{2}$ are prescribed, system \eqref{b} has important physical significance: in Bose-Einstein condensates and the nonlinear optics framework. When $V(x)=0$, $H_{u}(u,v)=\mu_{1}|u|^{p-2}u+\beta |u|^{\alpha-2} |v|^{\beta}u$ and $H_{v}(u,v)=\mu_{2}|v|^{q-2}v+\beta |u|^{\alpha}|v|^{\beta-2}v$, system \eqref{b} becomes the following type
\begin{equation}\label{c}
 \left\{
\begin{aligned}
&-\Delta u+\lambda_{1}u=\mu_{1}|u|^{p-2}u+\beta |u|^{\alpha-2} |v|^{\beta}u, \quad
\quad
\hbox{in }\mathbb{R}^{N},\\
&-\Delta v+\lambda_{2} v=\mu_{2}|v|^{q-2}v+\beta |u|^{\alpha}|v|^{\beta-2}v, \quad
\quad
\hbox{in }\mathbb{R}^{N},
\end{aligned}
\right.
\end{equation}
with the constraints
\begin{equation}\label{cc}
\int_{\mathbb{R}^{2}}|u|^{2}dx=a^{2}\quad \mbox{and} \quad \int_{\mathbb{R}^{2}}|v|^{2}dx=b^{2},
\end{equation}
where $\lambda_{1},\lambda_{2}$ cannot be prescribed but appear as Lagrange multipliers in a variational approach. For $N=3,4$, $p,q=2^{*}=\frac{2N}{N-2}$ and $\alpha,\beta>1$, Bartsch et al. \cite{BartschLi} proved existence and asymptotic behavior of normalized ground state solutions with $\alpha+\beta$ satisfy different conditions: $\alpha+\beta\leq 2+\frac{4}{N}$ or $\alpha+\beta> 2+\frac{4}{N}$. The most general form of \eqref{c} is $p=q=4$ and $\alpha=\beta=2$, which called the coupled Gross-Pitaevskii equations. Owing to it has a clear physical meaning, this case has been investigated by many authors. As is known to all, the results of system \eqref{c} with constraints \eqref{cc} have been widely studied in recent years; see e.g. \cite{BartschJeanjean,BartschJeanjeanSaove,BartschSaovea,BartschSaoveb,BartschZhong,Gou,Li} and references therein. However, in literature, there are no contribution devoted to the study of coupled systems with the constraints \eqref{cc} involving exponential nonlinearities.

When $V(x)=0$ and $u=v$, system \eqref{b} becomes the following type
\begin{equation}\label{equat1}
 \left\{
\begin{aligned}
&-\Delta u+\lambda u= f(u)  \quad \quad \hbox{in }\mathbb{R}^{N},\\
&u>0,\quad\quad\quad\int_{\mathbb{R}^{2}}|u|^{2}dx=a^{2}.
\end{aligned}
\right.
\end{equation}
For $N\geq1$, using the mountain pass lemma and a compactness argument, Jeanjean \cite{jeanjean1} studied the existence of normalized solutions for equation \eqref{equat1} with subcritical nonlinearities. In \cite{SOAVE1}, Soave considered the normalized solutions for equation \eqref{equat1} with combined subcritical nonlinearities where $f(t)=\mu|t|^{q-2}t+|t|^{p-2}t$ and $2<q\leq2+\frac{4}{N}\leq p <2^{*}$.

On the other hand, we mention that the study of our problem is based on some interesting results of critical case. For $N\geq3$ and $f$ has a critical growth in the Sobolev sense, Soave \cite{SOAVE} obtained the existence of ground states and raised some questions. For $N=2$, we observe $2^{*}=\infty$. Recently, Ji et al. \cite{AlvesJi} firstly studied equation \eqref{equat1} with exponential critical nonlinearities. They firstly obtained results regarding normalized problem with two-dimensional exponential critical growth. In addition, when $N\geq3$, the authors in \cite{AlvesJi} complemented some recent results found in \cite{SOAVE}.

Now, let us introduce the precise assumptions under what our problem is studied. Assume $H_{u},H_{v}$ are partial derivatives of a Carath\'{e}odory function $H$ and the nonlinearities $H_{u},H_{v}$ satisfy

\begin{itemize}
\item[\rm ($H_0$)]$H(w)\in C^{2}(\mathbb{R}^{2},\mathbb{R})$, where $w=(u,v)$;

\item[\rm ($H_1$)] $H_{u}(w)=o(|w|^{\tau})$ and $H_{v}(w)=o(|w|^{\tau})$ as $|w|\rightarrow0$, for some $\tau>3$;

\item[\rm ($H_2$)] there exists a positive constant $\theta>4$ such that $0<\theta H(w)\leq \nabla H(w)\cdot w$ for $a.e.$ $x\in \mathbb{R}^{2}$, where $\nabla H(w)=(H_{u}(w),H_{v}(w))$;

\item[\rm ($H_3$)]$H_{u}(u,0)=0$ for all $v\in \mathbb{R}$ and $H_{v}(0,v)=0$ for all $u\in \mathbb{R}$;

\item[\rm ($H_4$)] $H_{u}(w)u>0$ and $H_{v}(w)v>0$ for all $u,v\in \mathbb{R}$;

\item[\rm ($H_5$)]  Let $\tilde{H}(w)=\nabla H(w)\cdot w-2H(w)$ for all $w\in \mathbb{R}^{2}$. Then, $\nabla \tilde{H}(w)$ exists and
	\begin{equation*}
		\nabla \tilde{H}(w)\cdot w \geq 4 \tilde{H}(w), \ \mbox{for all} \ w \in \mathbb{R}^{2}.
	\end{equation*}

\item[\rm ($H_6$)]  there exist constants $\sigma>4$ and $\mu>0$ such that
	\begin{equation*}
		H(w)\geq \mu  \, |w|^{\sigma}\quad \text{for all}\,\, w\in \mathbb{R}^{2}.
    \end{equation*}

	\item[\rm ($H_7$)] $H_{u},H_{v}$ have $\gamma_{0}$-exponential critical growth, i.e. there exists $\gamma_{0}>0$ such that
\begin{align*}
    \lim\limits_{|w|\rightarrow +\infty}\frac{|H_{u}(w)|}{e^{\gamma|w|^{2}}}=\lim\limits_{|w|\rightarrow +\infty}\frac{|H_{v}(w)|}{e^{\gamma|w|^{2}}}=
 \left\{
    \begin{aligned}
    0, \ \ \ \forall \gamma> \gamma_{0},\\
    +\infty, \ \ \ \forall \gamma< \gamma_{0}.
    \end{aligned}
\right.
\end{align*}
\end{itemize}

Motivated by the research made in the critical Sobolev case and \cite{AlvesJi}, in this paper we consider the exponential critical growth for coupled system in $\mathbb{R}^{2}$. We recall that in $\mathbb{R}^{2}$, the natural growth restriction on the function $H_{u}$ and $H_{v}$ are given by the Trudinger-Moser inequality.

The result of this paper can be stated as follows:

\begin{theorem}\label{T2}
Assume that $H_{u},H_{v}$ satisfy ($H_{0}$)-($H_{4}$), $(H_6)$ and $(H_7)$. If $a^{2}+b^{2} <\frac{2\pi}{\gamma_{0}}$, then there exists $\mu_{1}=\mu_{1}(a,b)>0$ such that system \eqref{aa} has a weak solution $(\lambda_{1},\lambda_{2},u,v)$ with $\lambda_{1}>0$, $\lambda_{2}>0$ and $u,v\in H^{1}(\mathbb{R}^{2})$ for all $\mu \geq \mu_{1}$. Moreover, if $(H_{5})$ is also assumed, then $(u,v)$ can be chosen as a nontrivial ground state solution of system \eqref{aa}.
\end{theorem}

\begin{remark}\label{nontrivial}
If $H_{u},H_{v}$ dissatisfy ($H_{3}$) in Theorem \ref{T2}, we may obtain semitrivial ground state solutions of system \eqref{aa}. However, in this article, we require that $|u|_{2}=a>0$ and $|v|_{2}=b>0$. By $(H_{3})$, we know that $u=0$ implies that $v=0$. Thus, in Theorem \ref{T2} we obtain both $u$ and $v$ are nontrivial. In addition, in order to ensure that $|u|_{2}=a$ and $|v|_{2}=b$, the condition ($H_{4}$) is crucial.
\end{remark}

The paper is organized as follows: In Section $\ref{preliminaries}$, the variational setting and some preliminary results are presented. In Section $\ref{Sgeo}$, we introduce the geometry structure related to equation \eqref{aa}. In Section $\ref{PSa}$, we give some properties about the $(PS)$ sequence. Section $\ref{minimax}$ is devoted to give an estimate for the minimax level. Finally, in Section $\ref{critical}$, we prove Theorem \ref{T2}.

\vspace{0.5 cm}

\noindent \textbf{Notation:} From now on in this paper, otherwise mentioned, we use the following notations:
\begin{itemize}
    \item$C,C_1,C_2,...$ denote some positive constants.
	\item $|\,\,\,|_p$ denotes the usual norm of the Lebesgue space $L^{p}(\mathbb{R}^2)$, for $p \in [1,+\infty]$.
	\item $\Vert\,\,\,\Vert$ denotes the usual norm of the Sobolev space $H^{1}(\mathbb{R}^2)$.
    \item $o_{n}(1)$ denotes a real sequence with $o_{n}(1)\to 0$ as $n \rightarrow+\infty$.
\end{itemize}

\section{{\bfseries Preliminaries and functional setting}}\label{preliminaries}
In this section, we give some preliminary results and outline the variational framework for \eqref{aa}.
\begin{proposition}\label{PRa}
\cite{Cao}.$i)$ If $\gamma>0$ and $u\in H^{1}(\mathbb{R}^{2})$, then
\begin{eqnarray*}
    \begin{aligned}\displaystyle
    \int_{\mathbb{R}^{2}}\Big(e^{\gamma |u|^{2}}-1\Big)dx<\infty;
    \end{aligned}
\end{eqnarray*}
$ii)$ if $u\in H^{1}(\mathbb{R}^{2})$, $|\nabla u|^{2}_{2}\leq1$, $|u|_{2}\leq M<\infty$, and $\gamma<4\pi$, then there exists a constant $\mathcal{C}(M,\gamma)$, which depends only on $M$ and $\gamma$, such that
\begin{eqnarray*}
    \begin{aligned}\displaystyle
    \int_{\mathbb{R}^{2}}\Big(e^{\gamma |u|^{2}}-1\Big)dx<\mathcal{C}(M,\gamma).
    \end{aligned}
\end{eqnarray*}
\end{proposition}

Applying Proposition \ref{PRa}, we immediately have the following lemma.
\begin{lemma}\label{TMineq}
$i) $If $\gamma>0$ and $w=(u,v)\in H^{1}(\mathbb{R}^{2})\times H^{1}(\mathbb{R}^{2})$, then
\begin{eqnarray*}
    \begin{aligned}\displaystyle
    \int_{\mathbb{R}^{2}}\Big(e^{\gamma |w|^{2}}-1\Big)dx<\infty;
    \end{aligned}
\end{eqnarray*}
$ii)$ if $w=(u,v)\in H^{1}(\mathbb{R}^{2})\times H^{1}(\mathbb{R}^{2})$, $|\nabla w|^{2}_{2}\leq \frac{2\pi}{\gamma_{0}}-a^{2}-b^{2}$, $|u|_{2}^{2}= a^{2}$, $|v|_{2}^{2}= b^{2}$ and $\gamma<\gamma_{0}$ then there exists a constant $C( a,b,\gamma)$, which depends only on $a,b$ and $\gamma$ such that
\begin{eqnarray*}
    \begin{aligned}\displaystyle
    \int_{\mathbb{R}^{2}}\Big(e^{\gamma |w|^{2}}-1\Big)dx<C( a,b,\gamma).
    \end{aligned}
\end{eqnarray*}
\end{lemma}

\begin{proof}
By Proposition \ref{PRa} and Young's inequality, we obtain
\begin{eqnarray*}
    \begin{aligned}\displaystyle
    \int_{\mathbb{R}^{2}}\Big(e^{\gamma |w|^{2}}-1\Big)dx&=\int_{\mathbb{R}^{2}}\Big(e^{\gamma (u^{2}+v^{2})}-1\Big)dx\\
    &\leq \frac{1}{2} \int_{\mathbb{R}^{2}}\Big(e^{2\gamma u^{2}}-1\Big)dx+ \frac{1}{2} \int_{\mathbb{R}^{2}}\Big(e^{2\gamma v^{2}}-1\Big)dx\\
    &<\infty,
    \end{aligned}
\end{eqnarray*}
where $\gamma>0$ and $w=(u,v)\in H^{1}(\mathbb{R}^{2})\times H^{1}(\mathbb{R}^{2})$.

Moreover, by
\begin{eqnarray*}
    \begin{aligned}\displaystyle
	|\nabla w|_2^2 <\frac{2\pi}{\gamma_{0}}-a^{2}-b^{2} \quad \mbox{and} \quad |w|_{2}^{2}=|u|_{2}^{2}+|v|_{2}^{2}=a^2+b^{2},
    \end{aligned}
\end{eqnarray*}
we have
\begin{equation}\label{ener}
    \begin{aligned}\displaystyle
	2\gamma_{0}\|u\|^{2} < 4\pi \quad \mbox{and} \quad
	2\gamma_{0}\|v\|^{2} < 4\pi.
    \end{aligned}
\end{equation}
Then by Proposition \ref{PRa} and \eqref{ener}, we get
\begin{eqnarray*}
    \begin{aligned}\displaystyle
    \int_{\mathbb{R}^{2}}\Big(e^{\gamma |w|^{2}}-1\Big)dx&
    \leq \frac{1}{2} \int_{\mathbb{R}^{2}}\Big(e^{2\gamma  u^{2}}-1\Big)dx+ \frac{1}{2} \int_{\mathbb{R}^{2}}\Big(e^{2\gamma  v^{2}}-1\Big)dx\\
    &\leq \frac{1}{2} \int_{\mathbb{R}^{2}}\Big(e^{2\gamma_{0} \|u\|^{2} (\frac{|u|}{\|u\|})^{2}}-1\Big)dx +\frac{1}{2} \int_{\mathbb{R}^{2}}\Big(e^{2\gamma_{0} \|v\|^{2} (\frac{|v|}{\|v\|})^{2}}-1\Big)dx\\
    &\leq C,
    \end{aligned}
\end{eqnarray*}
where $C=C( a,b,\gamma)>0$. Then the proof is complete.
\end{proof}

Set $E=H^{1}(\mathbb{R}^{2})\times H^{1}(\mathbb{R}^{2})$, we define the following scalar product on $E$,
\begin{eqnarray*}
    \begin{aligned}\displaystyle
    \langle w_{1},w_{2}\rangle=\int_{\mathbb{R}^{2}}(\nabla u_{1}\nabla u_{2}+u_{1}u_{2})dx+\int_{\mathbb{R}^{2}}(\nabla v_{1}\nabla v_{2}+v_{1}v_{2})dx,
    \end{aligned}
\end{eqnarray*}
where $w_{1}=(u_{1},v_{1})$ and $w_{2}=(u_{2},v_{2})$, to which corresponds the norm
\begin{eqnarray*}
    \begin{aligned}\displaystyle
    \|w\|_{E}=\langle w,w\rangle^{1/2}.
    \end{aligned}
\end{eqnarray*}

Solutions of \eqref{aa} correspond to critical points of the energy functional $J:E\rightarrow \mathbb{R}$ defined by
\begin{equation}\label{func}
    \begin{aligned}\displaystyle
J(u,v)=\frac{1}{2}\int_{\mathbb{R}^{2}} (|\nabla u|^2 + |\nabla v|^2)dx-\int_{\mathbb{R}^{2}} H(u,v)dx
\end{aligned}
\end{equation}
and constrained to the $L^{2}$-torus
\begin{eqnarray*}
    \begin{aligned}\displaystyle
    \mathcal{T}(a,b):=\{(u,v)\in E: |u|_{2}=a, \ |v|_{2}=b\}.
    \end{aligned}
\end{eqnarray*}
The parameters $\lambda_{1},\lambda_{2}\in \mathbb{R}$ will appear as Lagrange multipliers. By Lemma \ref{TMineq}, we know that $H(w)\in L^{1}(\mathbb{R}^{2})$, which implies that $J$ is well defined. It is easy to see that $J$ is of class $C^{1}$, and that it is unbounded from below on $\mathcal{T}(a,b)$. It is well known that critical points of $J$ will not satisfy the Palais-Smale condition, as a consequence we recall that solutions of \eqref{aa} satisfy the Pohozaev identity
\begin{equation}\label{Pohozaev}
    \begin{aligned}\displaystyle
    P(u,v):=\int_{\mathbb{R}^{2}} (|\nabla u|^2 + |\nabla v|^2)dx-\int_{\mathbb{R}^{2}} \tilde{H}(u,v)dx,
    \end{aligned}
\end{equation}
where $\tilde{H}(u,v)=\nabla H(u,v)\cdot (u,v)-2H(u,v)$. Now, we introduce the $L^{2}$-invariant scaling $s\star u(x):=e^{s}u(e^{s}x)$ and $s\star (u,v):=(s\star u,s\star v)$. To be more precise, for $w=(u,v)\in \mathcal{T} (a,b)$ and $s\in \mathbb{R}$, let $\mathcal{F}: E\times \mathbb{R}\rightarrow E$ defined by
\begin{eqnarray*}
    \begin{aligned}\displaystyle
    \mathcal{F}(w, s)(x)=(e^{s}u(e^{s}x),e^{s}v(e^{s}x)).
    \end{aligned}
\end{eqnarray*}
Then $\tilde{J}:  \mathbb{R} \rightarrow \mathbb{R}$ with
\begin{eqnarray*}
    \begin{aligned}\displaystyle
	\tilde{J}_{w}(s):=J(\mathcal{F}(w, s))=\frac{e^{2s}}{2}\int_{\mathbb{R}^{2}} (|\nabla u|^2 +|\nabla v|^2 )dx-\frac{1}{e^{2s}}\int_{\mathbb{R}^{2}} H(e^{s}u(x),e^{s}v(x))\,dx.
    \end{aligned}
\end{eqnarray*}
A straightforward calculation shows that $\tilde{J}^{\prime}_{w}(s)=P(\mathcal{F}(w, s))$.

In order to overcome the loss of compactness of the Sobolev embedding in whole $\mathbb{R}^{2}$, we develop our variational procedure in the space $H^{1}_{rad}(\mathbb{R}^{2})$ to get some compactness results. Thus, we define
\begin{eqnarray*}
    \begin{aligned}\displaystyle
    E_{rad}:=H^{1}_{rad}(\mathbb{R}^{2})\times H^{1}_{rad}(\mathbb{R}^{2}), \ \ \mbox{and} \ \ \mathcal{T}_{r}(a,b):=\mathcal{T}(a,b)\cap E_{rad}.
    \end{aligned}
\end{eqnarray*}

From $(H_{1})$ and $(H_{7})$, fix $q>2$, for any $\varepsilon>0$ and $\gamma>\gamma_{0}$, there exist a constant $\kappa_{\varepsilon}>0$, which depends on $q,\gamma,\tau,\varepsilon$ and $\mu$, such that
\begin{equation}\label{nonb}
    \begin{aligned}\displaystyle
|\nabla H(w)|\leq|H_{u}(w)|+|H_{v}(w)|\leq\varepsilon|w|^{\tau}+\kappa_{\varepsilon}|w|^{q-1}(e^{\gamma |w|^{2}}-1) \ \mbox{for all} \ w\in \mathbb{R}^{2},
    \end{aligned}
\end{equation}
and by $(H_{2})$, we obtain
\begin{equation}\label{nona}
    \begin{aligned}\displaystyle
|H(w)|\leq\varepsilon|w|^{\tau+1}+\kappa_{\varepsilon}|w|^{q}(e^{\gamma |w|^{2}}-1) \ \mbox{for all} \ w\in \mathbb{R}^{2}.
    \end{aligned}
\end{equation}

\section{{\bfseries The minimax approach}}\label{Sgeo}To find a solution of \eqref{aa}, we show that $\tilde{J}$ on $\mathcal{T}_{r}(a,b)\times \mathbb{R}$ possesses a kind of mountain-pass geometrical structure.
\begin{lemma}\label{geo} Assume that $(H_{0})$-$(H_{2})$, $(H_{6})$ and $(H_{7})$ hold. Let $w\in \mathcal{T}_{r}(a,b)$. Then we have:

(i)$|\nabla\mathcal{F}(w, s)|_{2}\rightarrow 0$ and $J(\mathcal{F}(w, s))\rightarrow 0$ as $s\rightarrow -\infty$;

(ii)$|\nabla\mathcal{F}(w, s)|_{2}\rightarrow +\infty$ and $J(\mathcal{F}(w, s))\rightarrow -\infty$ as $s\rightarrow +\infty$.
\end{lemma}

\begin{proof} \mbox{}
 By a straightforward calculation, it follows that
\begin{equation}\label{cria}
    \int_{\mathbb{R}^2}|e^{s}u(e^{s}x)|^{2}\,dx=a^{2}, \quad \int_{\mathbb{R}^2}|e^{s}v(e^{s}x)|^{2}\,dx=b^{2},
\end{equation}
and
\begin{equation}\label{3.3}
    \int_{\mathbb{R}^2} |\nabla\mathcal{F}(w,s)(x)|^{2}\,dx=e^{2s}\int_{\mathbb{R}^2}(|\nabla u|^2+|\nabla v|^2 )dx.
\end{equation}
Moreover, we have
\begin{equation}\label{3.2}
\int_{\mathbb{R}^2}|\mathcal{F}(w, s)|^{\xi}\,dx= e^{(\xi-2)s}\int_{\mathbb{R}^2}|w(x)|^{\xi}\,dx, \quad \forall \xi \geq 2.
\end{equation}
From the above equalities, fixing $\xi>2$, we have
\begin{equation}\label{3.4}
    |\nabla\mathcal{F}(w,s)|_{2}\rightarrow 0 \quad \mbox{and} \quad |\mathcal{F}(w,s)|_{\xi} \rightarrow 0 \quad  \mbox{as} \quad s \to -\infty.
	\end{equation}
By \eqref{nona}, we obtain
\begin{eqnarray*}
    \begin{aligned}\displaystyle
|H(w)|\leq\varepsilon |w|^{\tau+1}+\kappa_{\varepsilon}|w|^{q}(e^{\gamma |w|^{2}}-1)\,\, \text{ for all }\, w \in E,
    \end{aligned}
\end{eqnarray*}
where $\gamma>\gamma_{0}$ close to $\gamma_{0}$ and $q>2$. Hence,
\begin{eqnarray*}
    \begin{aligned}\displaystyle
|H(\mathcal{F}(w,s))| \leq \varepsilon| \mathcal{F}(w, s)|^{\tau+1}+\kappa_{\varepsilon}| \mathcal{F}(w, s)|^{q}(e^{\gamma | \mathcal{F}(w, s)|^{2}}-1).
    \end{aligned}
\end{eqnarray*}
Using Lemma \ref{TMineq}, for $\gamma>\gamma_{0}$ close to $\gamma_{0}$, there exists $C=C(t, a,b,\gamma)>0$ such that
\begin{eqnarray*}
    \begin{aligned}\displaystyle
\int_{\mathbb{R}^2} (e^{\gamma | \mathcal{F}(w,s)|^{2}}-1)^{t}dx\leq\int_{\mathbb{R}^2} (e^{\gamma t | \mathcal{F}(w,s)|^{2}}-1)dx\leq C,
    \end{aligned}
\end{eqnarray*}
and so,
\begin{eqnarray*}
    \begin{aligned}\displaystyle
\int_{\mathbb{R}^2}|H(\mathcal{F}(w, s))|dx \leq \varepsilon\int_{\mathbb{R}^2}| \mathcal{F}(w, s)|^{\tau+1}dx+C_{1}\Big(\int_{\mathbb{R}^2}| \mathcal{F}(w,s)|^{qt'}dx\Big)^{1/t'},
    \end{aligned}
\end{eqnarray*}
where $t'=\frac{t}{t-1}$, and $t>1$ is close to $1$.  Now, by using (\ref{3.2}), we obtain
\begin{eqnarray*}
    \begin{aligned}\displaystyle
\int_{\mathbb{R}^2}|H(\mathcal{F}(w, s))|dx \leq \varepsilon e^{(\tau-1)s}\int_{\mathbb{R}^2} |w(x)|^{\tau+1}dx+C_{1}\Big(e^{(qt^{\prime}-2)s}\int_{\mathbb{R}^2}|w(x)|^{qt'}dx\Big)^{1/t'}.
    \end{aligned}
\end{eqnarray*}
Then, by $\tau-1>0$ and $qt^{\prime}-2>0$, we know that
\begin{eqnarray*}
    \begin{aligned}\displaystyle
	\int_{\mathbb{R}^2}|H( \mathcal{F}(w,s))| \rightarrow 0 \quad \mbox{as} \quad s \to -\infty,
    \end{aligned}
\end{eqnarray*}
which implies that $J(\mathcal{F}(w, s))\rightarrow 0$ as $s\rightarrow -\infty$, showing $(i)$.

In order to show $(ii)$, note that by (\ref{3.3}),
\begin{eqnarray*}
    \begin{aligned}\displaystyle
	|\nabla \mathcal{F}(w, s) |_{2}\rightarrow +\infty \quad \mbox{as} \quad s \to +\infty.
    \end{aligned}
\end{eqnarray*}
On the other hand, by $(H_{6})$, we obtain
\begin{eqnarray*}
    \begin{aligned}\displaystyle
	J(\mathcal{F}(w, s)) \leq \frac{1}{2}|\nabla \mathcal{F}(w,s)|_{2}^{2}-\frac{\mu}{\sigma}\int_{\mathbb{R}^2}|\mathcal{F}(w,s)|^{\sigma}\,dx=e^{2s}\int_{\mathbb{R}^2}|\nabla w|^2 dx-\frac{\mu e^{(\sigma-2)s}}{\sigma}\int_{\mathbb{R}^2}|w(x)|^{\sigma}\,dx.
    \end{aligned}
\end{eqnarray*}
Since $\sigma>4$, the last inequality yields $J(\mathcal{F}(w, s))\rightarrow -\infty$ as $s\rightarrow +\infty$.
\end{proof}

\begin{lemma}  \label{P1} Assume that $(H_{0})$-$(H_{2})$, $(H_{6})$ and $(H_{7})$ hold. Let $w\in \mathcal{T}_{r}(a,b)$. Then there exists $K(a,b)>0$ small enough such that
\begin{eqnarray*}
    \begin{aligned}\displaystyle
	0<\sup_{w\in A} J(w)<\inf_{w\in B} J(w)
    \end{aligned}
\end{eqnarray*}
	with
\begin{eqnarray*}
    \begin{aligned}\displaystyle
	A=\left\{w\in \mathcal{T}_{r}(a,b), \int_{\mathbb{R}^2} |\nabla w|^2 dx\leq K(a,b) \right\}\quad \mbox{and} \quad B=\left\{w\in \mathcal{T}_{r}(a,b), \int_{\mathbb{R}^2} |\nabla w|^2 dx=2K(a,b) \right\}.
    \end{aligned}
\end{eqnarray*}
\end{lemma}

\begin{proof}
	We will need the following Gagliardo-Sobolev inequality: for any $p> 2$,
\begin{eqnarray*}
    \begin{aligned}\displaystyle
	|u|_{p}\leq C(p, 2)|\nabla u|_2^{d_{p}}|u|_2^{1-d_{p}},
    \end{aligned}
\end{eqnarray*}
where $d_{p}=2(\frac{1}{2}-\frac{1}{p})$ and $u\in H^{1}(\mathbb{R}^{2})$. Then, let $w\in \mathcal{T}_{r}(a,b)$, we get an estimate of $|w|_{p}$. By Gagliardo-Sobolev inequality, there exist $C_{1},C_{2}>0$ such that
\begin{eqnarray*}
    \begin{aligned}\displaystyle
    |w|_{p}&=|u^{2}+v^{2}|_{\frac{p}{2}}^{\frac{1}{2}}\\
    &\leq (|u|_{p}^{2}+|v|_{p}^{2})^{\frac{1}{2}}\\
    &\leq (C_{1}|\nabla u|_{2}^{2d_{p}}|u|_{2}^{2-2d_{p}}+C_{2}|\nabla v|_{2}^{2d_{p}}|v|_{2}^{2-2d_{p}})^{\frac{1}{2}}\\
    &=(C_{1}a^{2-2d_{p}}|\nabla u|_{2}^{2d_{p}}+C_{2}b^{2-2d_{p}}|\nabla v|_{2}^{2d_{p}})^{\frac{1}{2}},
    \end{aligned}
\end{eqnarray*}
which implies that
\begin{equation}\label{GN}
    \begin{aligned}\displaystyle
|w|_{p}\leq C(|\nabla u|_{2}^{2d_{p}}+|\nabla v|_{2}^{2d_{p}})^{\frac{1}{2}}.
    \end{aligned}
\end{equation}
Now, let $K(a,b)<\frac{\pi}{\gamma_{0}}-\frac{a^{2}+b^{2}}{2}$. Assume that $w_{1},w_{2}\in \mathcal{T}_{r}(a,b)$ such that $|\nabla w_{1}|^{2}_{2}\leq K(a,b)$ and $|\nabla w_{2}|^{2}_2=2K(a,b)$. Thus, we have
\begin{eqnarray*}
    \begin{aligned}\displaystyle
\|u_2\|^{2}<\frac{2\pi}{\gamma_{0}} \quad \mbox{and} \quad  \|v_2\|^{2}<\frac{2\pi}{\gamma_{0}},
    \end{aligned}
\end{eqnarray*}
where $w_{2}=(u_2,v_2)$. Then, by Lemma \ref{TMineq}, we have
\begin{eqnarray*}
    \begin{aligned}\displaystyle
	\int_{\mathbb{R}^2}H(w_{2})\,dx \leq \varepsilon|w_{2}|^{\tau+1}_{\tau+1}+C_2|w_{2}|^{q}_{qt'},
    \end{aligned}
\end{eqnarray*}
where $\tau>3$, $q>2$, $t'=\frac{t}{t-1}>1$. Hence, by \eqref{GN}, we obtain
\begin{eqnarray*}
    \begin{aligned}\displaystyle
	\int_{\mathbb{R}^2}H(w_{2})\,dx \leq C_{1}K(a,b)^{\frac{\tau-1}{2}}+C_2K(a,b)^{(\frac{q}{2}-\frac{1}{t'})}.
    \end{aligned}
\end{eqnarray*}
From $(H_{2})$, we have $H(w_{1})> 0$ for any $w_{1}\in E$. Then
\begin{eqnarray*}
    \begin{aligned}\displaystyle
J(w_{2})-J(w_{1}) &=\frac{1}{2}\int_{\mathbb{R}^2}|\nabla w_{2}|^{2}\,dx-\frac{1}{2}\int_{\mathbb{R}^2}|\nabla w_{1}|^{2}\,dx-\int_{\mathbb{R}^2}H(w_{2})\,dx+\int_{\mathbb{R}^2}H(w_{1})\,dx\\
&\geq\frac{1}{2}\int_{\mathbb{R}^2}|\nabla w_{2}|^{2}\,dx-\frac{1}{2}\int_{\mathbb{R}^2}|\nabla w_{1}|^{2}\,dx-\int_{\mathbb{R}^2}H(w_{2})\,dx\\
&\geq K(a,b)-\frac{1}{2}K(a,b) -C_{3}K(a,b)^{\frac{\tau-1}{2}}-C_{4} K(a,b)^{(\frac{q}{2}-\frac{1}{t'})}.
    \end{aligned}
\end{eqnarray*}
Hence,
\begin{eqnarray*}
    \begin{aligned}\displaystyle
	J(w_{2})-J(w_{1}) \geq \frac{1}{2}K(a,b)-C_{3}K(a,b)^{\frac{\tau-1}{2}}-C_{4} K(a,b)^{(\frac{q}{2}-\frac{1}{t'})}.
    \end{aligned}
\end{eqnarray*}
Since $\tau>3$ and $t'>1$ with $\frac{q}{2}-\frac{1}{t'}>1$, choosing $K(a,b)$ small enough if necessary, it follows that
\begin{eqnarray*}
    \begin{aligned}\displaystyle
	\frac{1}{2}K(a,b)-C_3K(a,b)^{\frac{\tau-1}{2}}-C_4 K(a,b)^{(\frac{q}{2}-\frac{1}{t'})}>0,
    \end{aligned}
\end{eqnarray*}
	showing the desired result.
\end{proof}

As a byproduct of the last lemma is the following corollary.
\begin{corollary} \label{newcor}
Assume that $(H_{0})$-$(H_{2})$, $(H_{6})$ and $(H_{7})$ hold. Let $w\in \mathcal{T}_{r}(a,b)$. Then, if $|\nabla w|^{2}_{2}\leq K(a,b)$, there exists $K(a,b)>0$ small enough such that $J(w) >0$. Moreover,
\begin{eqnarray*}
    \begin{aligned}\displaystyle
	J_{\ast}=\inf\Big\{J(w): w\in \mathcal{T}_{r}(a,b) \ \ \mbox{and} \ \ |\nabla w|^{2}_{2}= \frac{K(a,b)}{2}\Big\}>0.
    \end{aligned}
\end{eqnarray*}
\end{corollary}

\begin{proof} Arguing as in the last lemma, we have for all $w\in \mathcal{T}_{r}(a,b)$
\begin{eqnarray*}
    \begin{aligned}\displaystyle
	J(w) \geq \frac{1}{2}|\nabla w|_{2}^{2}-C_{1}|\nabla w|_2^{\tau-1}-C_{2}|\nabla w|^{(q-\frac{2}{t^{\prime}})}_{2}>0,
    \end{aligned}
\end{eqnarray*}
	for $K(a,b)>0$ small enough.	
\end{proof}

In what follows, we fix $w_{0} \in \mathcal{T}_{r}(a,b)$ and apply Lemma \ref{geo}, \ref{P1} and Corollary  \ref{newcor}  to get two numbers $s_{1}<0$ and $s_{2}>0$, of such way that the functions $w_{1}=\mathcal{F}(w_{0}, s_{1})$ and $w_{2}=\mathcal{F}(w_{0}, s_{2})$ satisfy
\begin{eqnarray*}
    \begin{aligned}\displaystyle
|\nabla w_{1}|^2_2<\frac{K(a,b)}{2},\,\,  |\nabla w_{2}|_{2}^{2}>2K(a,b),\,\, J(w_{1})>0\,\,\mbox{and} \,\, J(w_{2})<0.
    \end{aligned}
\end{eqnarray*}

Now, following the ideas from Jeanjean \cite{jeanjean1}, we fix the following mountain pass level given by
\begin{eqnarray*}
    \begin{aligned}\displaystyle
m_\mu(a,b)=\inf_{h \in \Gamma}\max_{t \in [0,1]}J(h(t)),
    \end{aligned}
\end{eqnarray*}
where
\begin{eqnarray*}
    \begin{aligned}\displaystyle
\Gamma=\left\{h \in C([0,1], \ \mathcal{T}_{r}(a,b))\,:\,h(0)=w_{1} \,\, \mbox{and} \,\, h(1)=w_{2} \right\}.
    \end{aligned}
\end{eqnarray*}
By Corollary \ref{newcor}, we have
\begin{eqnarray*}
    \begin{aligned}\displaystyle
\max_{t \in [0,1]}J(h(t))\geq J_{\ast}>0.
    \end{aligned}
\end{eqnarray*}
Then, we obtain that $m_\mu(a,b)\geq J_{\ast}>0$.

\section{{\bfseries Palais-Smale sequence}}\label{PSa}
In this section, we take $\{w_{n}\}\subset \mathcal{T}_{r}(a,b)$ demotes the $(PS)$ sequence associated with the level $m_{\mu}(a,b)$ for $J$, where $w_{n}=(u_{n},v_{n})$. Using $\hat{w}_{n}=\mathcal{F}(w_{n},s_{n})$, we know that $\{\hat{w}_{n}\}$ is a $(PS)$ sequence associated with the level $m_{\mu}(a,b)$ for $\tilde{J}$. Thus, we have
\begin{equation} \label{mu2}
	J(w_{n}) \to m_\mu(a,b), \,\, \mbox{as} \,\, n \to +\infty,
\end{equation}
\begin{equation} \label{EQ2}
	-\Delta u_{n}+\lambda_{1,n}u_{n}=H_{u}(w_n)\ + o_n(1), \,\, \mbox{and} \ -\Delta v_{n}+\lambda_{2,n}v_{n}=H_{v}(w_n)\ + o_n(1),
\end{equation}
for some sequence $\{\lambda_{1,n}\},\{\lambda_{2,n}\} \subset \mathbb{R}$, and
\begin{equation} \label{PEQ2}
	P(w_{n})=\int_{\mathbb{R}^2}|\nabla w_{n}|^{2} dx +2 \int_{\mathbb{R}^2} H(w_{n}) dx- \int_{\mathbb{R}^2} \nabla H(w_{n}) \cdot w_{n} dx\to 0,\,\, \mbox{as} \,\, n \to +\infty.
\end{equation}
Then, setting
\begin{eqnarray*}
    \begin{aligned}\displaystyle
    \mathcal{P}(a,b):=\{(u,v)\in \mathcal{T}_{r}(a,b):P(u,v)=0\},
    \end{aligned}
\end{eqnarray*}
we are interested in the problem whether
\begin{eqnarray*}
    \begin{aligned}\displaystyle
    c(a,b):=\inf\limits_{(u,v)\in \mathcal{P}(a,b)}J(u,v)
    \end{aligned}
\end{eqnarray*}
is achieved. Then, in order to prove it, we give some properties of $(PS)$ sequence.

By \eqref{nona}, fixed $q>2$, for any $\varepsilon>0$ and $\gamma>\gamma_{0}$, there exist a constant $\kappa_{\varepsilon}>0$, which depends on $q,\gamma,\tau,\varepsilon$ and $\mu$, such that
\begin{equation}\label{nonlina}
|H(w)|\leq\varepsilon|w|^{\tau+1}+\kappa_{\varepsilon}|w|^{q}(e^{\gamma |w|^{2}}-1) \ \mbox{for all} \ w\in E.
\end{equation}
and
\begin{equation}\label{nonlinb}
|H(w_{n})|\leq\varepsilon|w_{n}|^{\tau+1}+\kappa_{\varepsilon}|w_{n}|^{q}(e^{\gamma |w_{n}|^{2}}-1), \ \mbox{for all} \ n\in \mathbb{N}.
\end{equation}

By \eqref{nonlina} and \eqref{nonlinb}, we obtain the convergence properties of the nonlinearities is related to the exponential function. Thus, the next lemma is important in our argument.

\begin{lemma}\label{imp} Let $\{w_{n}\}$ be a sequence in $E$ with $w_n=(u_{n},v_{n}) \in \mathcal{T}_{r}(a,b)$ and
\begin{equation}\label{energy}
    \begin{aligned}\displaystyle
	\limsup_{n \rightarrow +\infty} |\nabla w_n |_2^{2}  < \frac{2\pi}{\gamma_{0}}-a^{2}-b^{2}.
    \end{aligned}
\end{equation}
If $w_{n} \rightharpoonup w$ in $E_{rad}$ and $w_n(x) \rightarrow w(x)$ a.e. in $\mathbb{R}^{2}$, then
\begin{eqnarray*}
    \begin{aligned}\displaystyle
	|w_n|^{q}(e^{\gamma  |w_n(x)|^{2}}-1) \rightarrow |w|^{q}(e^{\gamma |w(x)|^{2}}-1) \,\, \mbox{in} \,\, L^{1}(\mathbb{R}^{2}).
    \end{aligned}
\end{eqnarray*}
\end{lemma}

\begin{proof}
Setting
\begin{eqnarray*}
    \begin{aligned}\displaystyle
	h_n(x)=e^{\gamma |w_n|^{2}}-1,
    \end{aligned}
\end{eqnarray*}
we can argue as in the proof of Lemma \ref{TMineq}, there exist $t>1$ and $t$ close to $1$  such that
\begin{eqnarray*}
    \begin{aligned}\displaystyle
    \int_{\mathbb{R}^{2}}\Big(e^{\gamma |w_{n}|^{2}}-1\Big)^{t}dx&
    \leq \int_{\mathbb{R}^{2}}\Big(e^{t\gamma |w_{n}|^{2}}-1\Big)dx\\&
    \leq \frac{1}{2} \int_{\mathbb{R}^{2}}\Big(e^{2t\gamma  u_{n}^{2}}-1\Big)dx+ \frac{1}{2} \int_{\mathbb{R}^{2}}\Big(e^{2t\gamma  v_{n}^{2}}-1\Big)dx\\
    &\leq \frac{1}{2} \int_{\mathbb{R}^{2}}\Big(e^{(2\gamma_{0}+\varepsilon)t \|u_{n}\|^{2} (\frac{|u_{n}|}{\|u_{n}\|})^{2}}-1\Big)dx +\frac{1}{2} \int_{\mathbb{R}^{2}}\Big(e^{(2\gamma_{0}+\varepsilon)t \|v_{n}\|^{2} (\frac{|v_{n}|}{\|v_{n}\|})^{2}}-1\Big)dx.
    \end{aligned}
\end{eqnarray*}
By \eqref{energy}, we have
\begin{eqnarray*}
    \begin{aligned}\displaystyle
	(2\gamma_{0}+\varepsilon)t\|u_n\|^{2} < 4\pi \quad \mbox{and} \quad
	(2\gamma_{0}+\varepsilon)t\|v_n\|^{2} < 4\pi.
    \end{aligned}
\end{eqnarray*}
Thus, there exists $C=C(t,a,b,\gamma)>0$ such that
\begin{eqnarray*}
    \begin{aligned}\displaystyle
   \int_{\mathbb{R}^{2}} \Big(h_n(x)\Big)^{t}dx=\int_{\mathbb{R}^{2}}\Big(e^{\gamma |w_{n}|^{2}}-1\Big)^{t}dx<C,
    \end{aligned}
\end{eqnarray*}
which implies that
\begin{eqnarray*}
    \begin{aligned}\displaystyle
	h_{n} \in L^{t}(\mathbb{R}^{2}) \quad \mbox{and} \quad \sup_{n \in \mathbb{N}}|h_n|_{t}<+\infty.
    \end{aligned}
\end{eqnarray*}
Therefore, $\{h_{n}\}$ is a bounded sequence in $L^{t}(\mathbb{R}^{2})$. Thus, for some subequence of $\{w_n\}$, still denoted by itself, we obtain that
	\begin{equation} \label{Newlimit}
		h_n \rightharpoonup h=e^{\gamma |w|^{2}}-1, \,\, \mbox{in} \,\, L^{t}(\mathbb{R}^{2}).
	\end{equation}
Now, we show that
\begin{equation}\label{ssubimp}
    \begin{aligned}\displaystyle
	|w_{n}|^{q}\rightarrow |w|^{q} \ \ \mbox{in} \ \ L^{t^{\prime}}(\mathbb{R}^{2}),
    \end{aligned}
\end{equation}
where $t^{\prime}=\frac{t}{t-1}$ and $w=(u,v)$. Then, by the embedding $H_{rad}^{1}(\mathbb{R}^{2}) \hookrightarrow L^{qt'}(\mathbb{R}^{2})$ is compact, we have
\begin{eqnarray*}
    \begin{aligned}\displaystyle
	u_{n} \rightarrow u \quad \mbox{and} \quad v_{n} \rightarrow v \quad \mbox{in} \quad L^{qt^{\prime}}(\mathbb{R}^{2}).
    \end{aligned}
\end{eqnarray*}
Thus, we obtain
\begin{equation}\label{sconveruv}
    \begin{aligned}\displaystyle
    |u_{n}|^{2} \rightarrow |u|^{2} \quad \mbox{and} \quad |v_{n}|^{2} \rightarrow |v|^{2} \quad \mbox{in} \quad L^{\frac{qt^{\prime}}{2}}(\mathbb{R}^{2}),
    \end{aligned}
\end{equation}
and
\begin{eqnarray*}
    \begin{aligned}\displaystyle
\Big||w_{n}|^{2}-|w|^{2}\Big|_{\frac{qt^{\prime}}{2}}&=\Big||u_{n}^{2}+v^{2}_{n}|-|u^{2}+v^{2}|\Big|_{\frac{qt^{\prime}}{2}}\\
&=\Big|(u_{n}^{2}-u^{2})+(v^{2}_{n}-v^{2})|\Big|_{\frac{qt^{\prime}}{2}}\\
&\leq |u_{n}^{2}-u^{2}|_{\frac{qt^{\prime}}{2}}+|v^{2}_{n}-v^{2}|_{\frac{qt^{\prime}}{2}}.
    \end{aligned}
\end{eqnarray*}
Then, using \eqref{sconveruv}, we know
\begin{eqnarray*}
    \begin{aligned}\displaystyle
 \Big||w_{n}|^{2}-|w|^{2}\Big|_{\frac{qt^{\prime}}{2}} \leq |u_{n}^{2}-u^{2}|_{\frac{qt^{\prime}}{2}}+|v^{2}_{n}-v^{2}|_{\frac{qt^{\prime}}{2}} \ \rightarrow 0 \ \mbox{as} \ n\rightarrow\infty,
    \end{aligned}
\end{eqnarray*}
which implies that
\begin{eqnarray*}
    \begin{aligned}\displaystyle
    |w_{n}|^{2}\rightarrow |w|^{2} \ \ \mbox{in} \ \ L^{\frac{qt^{\prime}}{2}}(\mathbb{R}^{2}).
    \end{aligned}
\end{eqnarray*}
Hence,
\begin{equation}\label{subimp}
    \begin{aligned}\displaystyle
	|w_{n}|^{q}\rightarrow |w|^{q} \ \ \mbox{in} \ \ L^{t^{\prime}}(\mathbb{R}^{2}),
    \end{aligned}
\end{equation}
where $t^{\prime}=\frac{t}{t-1}$ and $w=(u,v)$. Together \eqref{Newlimit} with \eqref{subimp}, we know
\begin{eqnarray*}
    \begin{aligned}\displaystyle
    |w_n|^{q}(e^{\gamma  |w_n(x)|^{2}}-1) \rightarrow |w|^{q}(e^{\gamma |w(x)|^{2}}-1) \,\, \mbox{in} \,\, L^{1}(\mathbb{R}^{2}).
    \end{aligned}
\end{eqnarray*}
Then, the proof is complete.
\end{proof}

By Lemma \ref{imp}, we have the following two important Corollaries.
\begin{corollary}\label{slimitb}Assume that $(H_{0})$-$(H_{2})$, $(H_{6})$ and $(H_{7})$ hold. Let $\{w_{n}\}$ be the $(PS)$ sequence of $J$ with $w_{n}=(u_{n},v_{n}) \in \mathcal{T}_{r}(a,b)$ and
\begin{eqnarray*}
    \begin{aligned}\displaystyle
	\limsup_{n \rightarrow +\infty} |\nabla w_{n} |_2^{2}  < \frac{2\pi}{\gamma_{0}}-a^2-b^{2}.
    \end{aligned}
\end{eqnarray*}
	If $w_{n} \rightharpoonup w$ in $E_{rad}$ and $w_{n}(x) \rightarrow w(x)$ a.e in $\mathbb{R}^{2}$, then, for any $\phi,\psi\in C_{0}^{\infty}(\mathbb{R}^{2})$, we have
\begin{eqnarray*}
    \begin{aligned}\displaystyle
H_{u}(w_n)\phi \to H_{u}(w)\phi \quad \mbox{and} \quad H_{v}(w_n)\psi \to H_{v}(w)\psi \,\, \mbox{in} \,\, L^{1}(\mathbb{R}^2).
    \end{aligned}
\end{eqnarray*}
\end{corollary}

\begin{proof}
By \eqref{nonb}, we know that
\begin{eqnarray*}
    \begin{aligned}\displaystyle
|H_{u}(w_{n})|\leq\varepsilon|w_{n}|^{\tau}+\kappa_{\varepsilon}|w_{n}|^{q-1}(e^{\gamma |w_{n}|^{2}}-1) \ \mbox{for all} \ w_{n}\in E.
    \end{aligned}
\end{eqnarray*}
Hence, for any $\phi\in C_{0}^{\infty}(\mathbb{R}^{2})$, we have
\begin{eqnarray*}
    \begin{aligned}\displaystyle
|H_{u}(w_{n})\phi|\leq |H_{u}(w_{n})||\phi|\leq   \varepsilon|w_{n}|^{\tau}|\phi|+\kappa_{\varepsilon}|w_{n}|^{q-1}|\phi|(e^{\gamma |w_{n}|^{2}}-1) \ \mbox{for all} \ w_{n}\in E.
    \end{aligned}
\end{eqnarray*}
Let $U=supp \phi$. Then, we obtain
\begin{eqnarray*}
    \begin{aligned}\displaystyle
    \int_{U}|u_{n}|^{\tau}|\phi|dx\rightarrow\int_{U}|u|^{\tau}|\phi|dx, \ \mbox{as} \ n\rightarrow\infty,
    \end{aligned}
\end{eqnarray*}
and
\begin{eqnarray*}
    \begin{aligned}\displaystyle
    \int_{U}|u_{n}|^{q-1}|\phi|(e^{\gamma |u_{n}|^{2}}-1)dx\rightarrow\int_{U}|u|^{q-1}|\phi|(e^{\gamma |u|^{2}}-1)dx, \ \mbox{as} \ n\rightarrow\infty.
    \end{aligned}
\end{eqnarray*}
Now, applying a variant of the Lebesgue Dominated Convergence Theorem, we can deduce that
\begin{eqnarray*}
    \begin{aligned}\displaystyle
H_{u}(w_n)\phi \to H_{u}(w)\phi \,\, \mbox{in} \,\, L^{1}(\mathbb{R}^2), \ \mbox{for any }\phi\in C_{0}^{\infty}(\mathbb{R}^{2}).
    \end{aligned}
\end{eqnarray*}
A similar argument works to show that
\begin{eqnarray*}
    \begin{aligned}\displaystyle
H_{v}(w_n)\psi \to H_{v}(w)\psi \,\, \mbox{in} \,\, L^{1}(\mathbb{R}^2), \ \mbox{for any }\psi\in C_{0}^{\infty}(\mathbb{R}^{2}),
    \end{aligned}
\end{eqnarray*}
which completes the proof.
\end{proof}

\begin{corollary} \label{Conver}Assume that $(H_{0})$-$(H_{2})$, $(H_{6})$ and $(H_{7})$ hold. Let $\{w_{n}\}$ be the $(PS)$ sequence of $J$ with $w_{n}=(u_{n},v_{n}) \in \mathcal{T}_{r}(a,b)$ and
\begin{eqnarray*}
    \begin{aligned}\displaystyle
	\limsup_{n \rightarrow +\infty} |\nabla w_{n} |_2^{2}  < \frac{2\pi}{\gamma_{0}}-a^2-b^{2}.
    \end{aligned}
\end{eqnarray*}
	If $w_{n} \rightharpoonup w$ in $E_{rad}$ and $w_{n}(x) \rightarrow w(x)$ a.e in $\mathbb{R}^{2}$, then
\begin{eqnarray*}
    \begin{aligned}\displaystyle
	H(w_n) \to H(w) \quad \mbox{and} \quad \nabla H(w_{n})\cdot w_{n} \to \nabla H(w)\cdot w \,\, \mbox{in} \,\, L^{1}(\mathbb{R}^2).
    \end{aligned}
\end{eqnarray*}	
\end{corollary}
\begin{proof}By \eqref{nona}, we have
\begin{eqnarray*}
    \begin{aligned}\displaystyle
	|H(w)|\leq\varepsilon|w|^{\tau+1}+C|w|^{q}(e^{\gamma |w|^{2}}-1)\,\, \text{ for all }\, w \in E,
    \end{aligned}
\end{eqnarray*}
where $\gamma>\gamma_{0}$, $\tau>3$ and $q>2$. Therefore,
\begin{equation} \label{Domina1}
		|H(w_{n})|\leq\varepsilon|w_{n}|^{\tau+1}+C|w_{n}|^{q}(e^{\gamma |w_{n}|^{2}}-1)\,\, \text{ for all }\, w_{n} \in E.
	\end{equation}
By 	Lemma \ref{imp}, we know
\begin{eqnarray*}
    \begin{aligned}\displaystyle
	|w_{n}|^{q}(e^{\gamma |w_{n}(x)|^{2}}-1) \rightarrow |w|^{q}(e^{\gamma |w(x)|^{2}}-1) \quad \mbox{in} \quad L^{1}(\mathbb{R}^{2}).
    \end{aligned}
\end{eqnarray*}
By the compact embedding $H_{rad}^{1}(\mathbb{R}^{2}) \hookrightarrow L^{p}(\mathbb{R}^{2})$, we have
\begin{eqnarray*}
    \begin{aligned}\displaystyle
	u_{n}\rightarrow u \quad \mbox{and} \quad v_{n}\rightarrow v \quad \mbox{in} \ \ L^{p}(\mathbb{R}^2).
    \end{aligned}
\end{eqnarray*}
Now, by \eqref{nonlina} and \eqref{nonlinb}, we can use  a variant of the Lebesgue Dominated Convergence Theorem to conclude that
\begin{eqnarray*}
    \begin{aligned}\displaystyle
	H(w_{n}) \to H(w) \quad \mbox{in} \quad L^{1}(\mathbb{R}^2).
    \end{aligned}
\end{eqnarray*}
	A similar argument works to show that
\begin{eqnarray*}
    \begin{aligned}\displaystyle
	\nabla H(w_{n})\cdot w_{n} \to \nabla H(w)\cdot w \quad \mbox{in} \quad L^{1}(\mathbb{R}^2),
    \end{aligned}
\end{eqnarray*}	
which completes the proof.
\end{proof}

\begin{lemma}\label{uni}
Assume that ($H_{0}$) and ($H_{5}$) hold. Let $w\in \mathcal{T}_{r}(a,b)$. Then the function $\tilde{J}_{w}(s)=J(\mathcal{F}(w,s))$ reaches its unique maximum at a point $s(w)\in \mathbb{R}$ such that $\mathcal{F}(w,s(w))\in \mathcal{P}(a,b)$.
\end{lemma}
\begin{proof}
First, we know
\begin{equation}\label{DS1}
    \begin{aligned}\displaystyle
\tilde{J}^{\prime}_{w}(s)=P(\mathcal{F}(w,s))=|\nabla \mathcal{F}(w,s)|^{2}_{2}-\int_{\mathbb{R}^{2}}\tilde{H}(\mathcal{F}(w,s)(x))dx.
    \end{aligned}
\end{equation}

By Lemmas \ref{geo} and \ref{P1}, we know that there exists at least a $s_{0}\in \mathbb{R}$ such that $\tilde{J}^{\prime}_{w}(s)_{|_{s=s_{0}}}=0$. By \eqref{DS1}, we obtain that $\mathcal{F}(w,s_{0})\in \mathcal{P}(a,b)$. Then, by $\tilde{J}^{\prime}_{w}(s)_{|_{s=s_{0}}}=0$, we have
\begin{eqnarray*}
    \begin{aligned}\displaystyle
\tilde{J}^{\prime\prime}_{w}(s)_{|_{s=s_{0}}}= 4\int_{\mathbb{R}^{2}}\tilde{H}(\mathcal{F}(w,s_{0})(x))dx -\int_{\mathbb{R}^{2}}\nabla \tilde{H}(\mathcal{F}(w,s_{0})(x))\cdot \mathcal{F}(w,s_{0})(x)dx.
    \end{aligned}
\end{eqnarray*}
Thus, by ($H_{5}$), we deduce that $\tilde{J}^{\prime\prime}_{w}(s)_{|_{s=s_{0}}}<0$, which proves the unicity of $s_{0}$.
\end{proof}

By virtue of Lemma \ref{uni},  similarly to \cite[Lemma 2.10]{jeanjean1} or \cite[Lemma 2.5]{BartschLi}, we show the following Lemma.

\begin{lemma}\label{fin}
Assume that ($H_{0}$) and ($H_{5}$) hold. Let $w\in \mathcal{T}_{r}(a,b)$. Then
\begin{eqnarray*}
    \begin{aligned}\displaystyle
    m_{\mu}(a,b)= c(a,b):=\inf\limits_{w\in \mathcal{P}(a,b)}\max\limits_{s\in \mathbb{R}}J(\mathcal{F}(w,s)).
    \end{aligned}
\end{eqnarray*}
\end{lemma}

\begin{proof}
For $w\in \mathcal{T}_{r}(a,b)$, by Lemma \ref{uni}, we know that
\begin{eqnarray*}
    \begin{aligned}\displaystyle
    \max\limits_{s\in \mathbb{R}}J(\mathcal{F}(w,s))=J(\mathcal{F}(w,s(w)))\geq m_{\mu}(a,b),
    \end{aligned}
\end{eqnarray*}
Thus, $m_{\mu}(a,b)\leq c(a,b)$. On the other hand, for $w\in \mathcal{P}(a,b)$, Lemma \ref{uni} also implies
\begin{eqnarray*}
    \begin{aligned}\displaystyle
J(w)=\max\limits_{s\in \mathbb{R}}J(\mathcal{F}(w,s))\geq \inf\limits_{w\in \mathcal{P}(a,b)}\max\limits_{s\in \mathbb{R}}J(\mathcal{F}(w,s)),
    \end{aligned}
\end{eqnarray*}
hence, $m_{\mu}(a,b)= c(a,b)$.
\end{proof}

\section{{\bfseries On the mini-max level}}\label{minimax}
In this section, we obtain an upper bound for the minimax level. Thus, we obtain an upper bound for $\|w\|_{E}^{2}$, which is important for exponential critical problem.
\begin{lemma} \label{ESTMOUNTPASS} Assume that ($H_{6}$) holds, there holds $\displaystyle \lim_{\mu \rightarrow +\infty}m_\mu(a,b)=0$.	
\end{lemma}

\begin{proof} In what follow we set the path $h_{0}(t)=\mathcal{F}\big(w_{0}, (1-t)s_1+ts_2\big) \in \Gamma$. Then, by $(H_{7})$,
\begin{eqnarray*}
    \begin{aligned}\displaystyle
	m_\mu(a,b) \leq \max_{t \in [0,1]}J(h_0(t))&\leq \max_{s \in \mathbb{R}}\left\{\frac{e^{2s}}{2}|\nabla w_{0}|_{2}^{2}-\frac{\mu}{\sigma}e^{(\sigma-2)s}|w_{0}|_{\sigma}^{\sigma}\right\}\\
 &= \max_{r \geq 0}\left\{\frac{r^{2}}{2}|\nabla w_{0}|_{2}^{2}-\frac{\mu}{\sigma}r^{\sigma-2}|w_{0}|_{\sigma}^{\sigma}\right\}.
    \end{aligned}
\end{eqnarray*}
Thus,
\begin{eqnarray*}
    \begin{aligned}\displaystyle
	m_\mu(a,b) \leq C_{2}\left(\frac{1}{\mu}\right)^{\frac{2}{\sigma-4}} \to 0 \quad \mbox{as} \quad \mu \to +\infty,
    \end{aligned}
\end{eqnarray*}
for some $C_2>0$. Here, we have used the fact that $\sigma>4$.
\end{proof}

\begin{lemma}\label{newlem1} Assume that $(H_{2})$ and ($H_{6}$) hold. Let $w_{n}\in \mathcal{T}_{r}(a,b)$. Then, the $(PS)$ sequence $\{w_{n}\}$ of $J$ satisfies
\begin{eqnarray*}
    \begin{aligned}\displaystyle
	\limsup_{n \rightarrow +\infty}\int_{\mathbb{R}^2}H(w_{n})\,dx \leq \frac{2}{\theta-4}m_\mu(a,b).
    \end{aligned}
\end{eqnarray*}
\end{lemma}
\begin{proof} Using the fact that $J(w_{n})=m_\mu(a,b)+o_n(1)$ and $P(w_{n})=o_n(1)$, it follows that
\begin{eqnarray*}
    \begin{aligned}\displaystyle
	2{J}(w_{n})+P(w_{n})=2m_\mu(a,b)+o_n(1),
    \end{aligned}
\end{eqnarray*}
and so,
\begin{equation}\label{equa1}
    \begin{aligned}\displaystyle
	2|\nabla w_{n}|_{2}^{2}-\int_{\mathbb{R}^2}\nabla H(w_{n})\cdot w_{n}\,dx=2m_\mu(a,b)+o_n(1).
    \end{aligned}
\end{equation}
Using that $J(w_{n})=m_\mu(a,b)+o_n(1)$, we obtain
\begin{equation}\label{equa2}
    \begin{aligned}\displaystyle
|\nabla w_{n}|_{2}^{2}=2\int_{\mathbb{R}^2} H(w_{n}) dx+2m_\mu(a,b)+o_n(1).
    \end{aligned}
\end{equation}
Together \eqref{equa1} and \eqref{equa2}, we get
\begin{eqnarray*}
    \begin{aligned}\displaystyle
4\int_{\mathbb{R}^2} H(w_{n}) dx+4m_\mu(a,b)+o_n(1)-\int_{\mathbb{R}^2}\nabla H(w_{n})\cdot w_{n}\,dx=2m_\mu(a,b)+o_n(1).
    \end{aligned}
\end{eqnarray*}
Then, by $(H_{2})$,
\begin{eqnarray*}
    \begin{aligned}\displaystyle
2m_\mu(a,b)+o_n(1)=\int_{\mathbb{R}^2}\nabla H(w_{n})\cdot w_{n}\,dx-4\int_{\mathbb{R}^2} H(w_{n}) dx\geq (\theta-4)\int_{\mathbb{R}^2} H(w_{n}) dx.
    \end{aligned}
\end{eqnarray*}
Since $\theta>4$, we have
\begin{eqnarray*}
    \begin{aligned}\displaystyle
\limsup_{n \rightarrow +\infty}\int_{\mathbb{R}^2}H(w_{n})\,dx \leq \frac{2}{\theta-4}m_\mu(a,b).
    \end{aligned}
\end{eqnarray*}
Then, the proof is complete.
\end{proof}

\begin{lemma} \label{boundd}Assume that $(H_{2})$ and ($H_{6}$) hold. Let $w_{n}\in \mathcal{T}_{r}(a,b)$. Then, the $(PS)$ sequence $\{w_{n}\}$ of $J$ satisfies
\begin{eqnarray*}
    \begin{aligned}\displaystyle
\displaystyle \limsup_{n \rightarrow +\infty}|\nabla w_{n}|_{2}^{2} \leq \frac{2(\theta-2)}{\theta-4} m_\mu(a,b).
    \end{aligned}
\end{eqnarray*}
Hence, there exists $\mu^*>0$ such that
\begin{equation}\label{level}
    \begin{aligned}\displaystyle
	\limsup_{n \to +\infty}|\nabla w_{n}|_{2}^{2}<\frac{2\pi}{\gamma_{0}}-a^{2}-b^{2}, \,\,\text{for any} \,\,\,\mu \geq \mu^*.
    \end{aligned}
\end{equation}
Moreover, the sequence $\{w_{n}\}$ is bounded in $E$.
\end{lemma}

\begin{proof}  Since $J(w_{n})=m_\mu(a,b)+o_n(1)$, we have
\begin{eqnarray*}
    \begin{aligned}\displaystyle
	\int_{\mathbb{R}^2}|\nabla w_{n}|^2\,dx=2m_\mu(a,b)+2\int_{\mathbb{R}^2}H(w_n)\,dx+o_n(1).
    \end{aligned}
\end{eqnarray*}
Thus, by Lemma \ref{newlem1}, we obtain
\begin{eqnarray*}
    \begin{aligned}\displaystyle
	\limsup_{n \to +\infty}|\nabla w_n|_{2}^{2}  \leq \frac{2(\theta-2)}{\theta-4} m_\mu(a,b),
    \end{aligned}
\end{eqnarray*}
which implies \eqref{level} hold and $\{w_{n}\}$ is bounded in $E$.
\end{proof}

\begin{lemma} \label{mu}Assume that $(H_{2})$, $(H_{4})$ and ($H_{6}$) hold. For all $\mu \geq \mu_{1}$, where $\mu_{1}=\max\{\mu^{*},\mu^{**}\}$, we have  $\{\lambda_{1,n}\},\{\lambda_{2,n}\}$ are bounded sequence with
\begin{eqnarray*}
    \begin{aligned}\displaystyle
	 \lambda_{1,n}\rightarrow \lambda_{1}>0,  \lambda_{2,n}\rightarrow \lambda_{2}>0 \ \mbox{and} \ \limsup_{n \rightarrow +\infty} |-\lambda_{1,n}a^{2}-\lambda_{2,n}b^{2}|\leq \frac{4(\theta-1)}{\theta-4} m_\mu(a,b).
    \end{aligned}
\end{eqnarray*}
\end{lemma}

\begin{proof}
By \eqref{EQ2}, the number $\lambda_{1,n},\lambda_{2,n}$  must satisfy the equality below
\begin{equation} \label{lambdan1}
	\lambda_{1,n}=-\frac{1}{a^{2}}\left\{ J^{\prime}(u_{n},v_{n})(u_{n},0) \right\}+o_n(1) \ \mbox{and} \ \lambda_{2,n}=-\frac{1}{b^{2}}\left\{ J^{\prime}(u_{n},v_{n})(0,v_{n}) \right\}+o_n(1),
\end{equation}
thus, $\{\lambda_{1,n}\},\{\lambda_{2,n}\}$ are bounded in $\mathbb{R}$. Then, by \eqref{lambdan1}, we obtain
\begin{equation}\label{lama}
    \begin{aligned}\displaystyle
    -\lambda_{1,n}a^{2} =|\nabla u_{n}|_{2}^{2}-\int_{\mathbb{R}^2} H_{u}(w_{n}) u_{n}\,dx+o_n(1),
    \end{aligned}
\end{equation}
and
\begin{equation}\label{lamb}
    \begin{aligned}\displaystyle
    -\lambda_{2,n}b^{2} =|\nabla v_{n}|_{2}^{2}-\int_{\mathbb{R}^2} H_{v}(w_{n}) v_{n}\,dx+o_n(1).
    \end{aligned}
\end{equation}
The limit (\ref{PEQ2}) together with Lemmas \ref{newlem1} and \ref{boundd} ensures that $(\int_{\mathbb{R}^2}\nabla H(w_{n})\cdot w_{n}\,dx)$ is bounded with
\begin{equation}\label{dnon}
    \begin{aligned}\displaystyle
\limsup_{n \rightarrow +\infty}\int_{\mathbb{R}^2}\nabla H(w_{n})\cdot w_{n}\,dx\leq \frac{2\theta}{\theta-4} m_\mu(a,b).
    \end{aligned}
\end{equation}
By \eqref{lama} and Lemma \ref{boundd}, we obtain
\begin{eqnarray*}
    \begin{aligned}\displaystyle
    \lambda_{1,n}a^{2}&=-|\nabla u_{n}|_{2}^{2}+\int_{\mathbb{R}^2} H_{u}(w_{n}) u_{n}\,dx+o_n(1)\\
    &\geq-|\nabla w_{n}|_{2}^{2} +\int_{\mathbb{R}^2} H_{u}(w_{n}) u_{n}\,dx+o_n(1)\\
    &\geq-\frac{2(\theta-2)}{\theta-4} m_\mu(a,b)+\int_{\mathbb{R}^2} H_{u}(w_{n}) u_{n}\,dx+o_n(1).
    \end{aligned}
\end{eqnarray*}
Then, by $(H_{4})$, there exists $\mu^{**}>0$ such that
\begin{equation}\label{lamdda}
    \begin{aligned}\displaystyle
\lambda_{1,n}a^{2}\geq-\frac{2(\theta-2)}{\theta-4} m_\mu(a,b)+\int_{\mathbb{R}^2} H_{u}(w_{n}) u_{n}\,dx+o_n(1)>0 \ \mbox{for all} \  \mu \geq \mu^{**}.
    \end{aligned}
\end{equation}
A similar argument works to show that
\begin{equation}\label{lamddb}
    \begin{aligned}\displaystyle
\lambda_{2,n}b^{2}\geq-\frac{2(\theta-2)}{\theta-4} m_\mu(a,b)+\int_{\mathbb{R}^2} H_{v}(w_{n}) v_{n}\,dx+o_n(1)>0 \ \mbox{for all} \  \mu \geq \mu^{**}.
    \end{aligned}
\end{equation}
By $\{\lambda_{1,n}\},\{\lambda_{2,n}\}$ are bounded in $\mathbb{R}$, using \eqref{lamdda} and \eqref{lamddb}, we have
\begin{eqnarray*}
    \begin{aligned}\displaystyle
    \lambda_{1,n}\rightarrow \lambda_{1}>0, \ \mbox{and} \ \lambda_{2,n}\rightarrow \lambda_{2}>0.
    \end{aligned}
\end{eqnarray*}

Then, combining the two equality of \eqref{lambdan1}, we obtain
\begin{eqnarray*}
    \begin{aligned}\displaystyle
	-\lambda_{1,n}a^{2}-\lambda_{2,n}b^{2} =|\nabla u_{n}|_{2}^{2}+|\nabla v_{n}|_{2}^{2}-\int_{\mathbb{R}^2}\nabla H(u_{n},v_{n})\cdot (u_{n},v_{n})\,dx+o_n(1).
    \end{aligned}
\end{eqnarray*}
Hence, we obtain
\begin{eqnarray*}
    \begin{aligned}\displaystyle
\limsup_{n \rightarrow +\infty} |-\lambda_{1,n}a^{2}-\lambda_{2,n}b^{2}|&\leq
|\nabla w_{n}|_{2}^{2}+\int_{\mathbb{R}^2}\nabla H(w_{n})\cdot w_{n}\,dx\\
&\leq \frac{2\theta}{(\theta-4)}m_\mu(a,b)+\frac{2(\theta-2)}{\theta-4} m_\mu(a,b)\\&=\frac{4(\theta-1)}{\theta-4} m_\mu(a,b),
    \end{aligned}
\end{eqnarray*}
which completes the proof.
\end{proof}

\section{{\bfseries Proof of Theorem \ref{T2}}}\label{critical}
In this section, we restrict our study to the space $E_{rad}$ and assume that $H_{u},H_{v}$ have critical growth with exponent critical $\gamma_{0}$.\\
{\bf Proof of Theorem}$\ref{T2}$. First, by Lemma \ref{mu}, we have $\lambda_{1}>0$ and $\lambda_{2}>0$. Then, the equality \eqref{EQ2} and Corollary \ref{slimitb} imply that
\begin{eqnarray*}
    \begin{aligned}\displaystyle
-\Delta u+\lambda_{1}u=H_{u}(w) \ \ \mbox{and} \ -\Delta v+\lambda_{1}v=H_{v}(w) \ \  \mbox{in} \ H^{1}(\mathbb{R}^{2}).
    \end{aligned}
\end{eqnarray*}
Thus, $P(u,v)=0$. Then, we show $w_{n}\rightharpoonup w:=(u,v)$ in $E_{rad}$, where $u,v\neq 0$. By Lemma \ref{boundd}, we have $w_{n}\rightharpoonup w:=(u,v)$ in $E_{rad}$. Assume that $u=0$, then by ($H_{3}$) and Remark \ref{nontrivial}, we know that $v=0$. Thus, the only case is that $w_{n}\rightharpoonup w=(0,0)$. For any $\mu \geq \mu_{1}$, using Corollary \ref{Conver}, it follows that
\begin{eqnarray*}
    \begin{aligned}\displaystyle
\lim_{n \rightarrow +\infty}\int_{\mathbb{R}^2}\nabla H(w_{n})\cdot w_{n}\,dx=\int_{\mathbb{R}^2}\nabla H(w)\cdot w\,dx,
    \end{aligned}
\end{eqnarray*}
and
\begin{eqnarray*}
    \begin{aligned}\displaystyle
\lim_{n \rightarrow +\infty}\int_{\mathbb{R}^2}H(w_{n})\,dx=\int_{\mathbb{R}^2}H(w)\,dx,
    \end{aligned}
\end{eqnarray*}
where $w_{n} \rightharpoonup w$ in $E_{rad}$. The last limit implies that $w \not=(0,0)$, because otherwise, Corollary \ref{Conver} gives
\begin{eqnarray*}
    \begin{aligned}\displaystyle
\lim_{n \rightarrow +\infty}\int_{\mathbb{R}^2}H(w_{n})\,dx=\lim_{n \rightarrow +\infty}\int_{\mathbb{R}^2}\nabla H(w_{n})\cdot w_{n}\,dx=0,
    \end{aligned}
\end{eqnarray*}
and by Lemma \ref{mu},
\begin{eqnarray*}
    \begin{aligned}\displaystyle
\limsup_{n \rightarrow +\infty} |-\lambda_{1,n}a^{2}-\lambda_{2,n}b^{2}|\leq \frac{4(\theta-1)}{\theta-4} m_\mu(a,b),
    \end{aligned}
\end{eqnarray*}
which implies that
\begin{eqnarray*}
    \begin{aligned}\displaystyle
\limsup_{n \rightarrow +\infty} (-\lambda_{1,n}a^{2}-\lambda_{2,n}b^{2})\leq0.
    \end{aligned}
\end{eqnarray*}
Since $\{w_{n}\}$ is bounded in $E$ and $\displaystyle \limsup_{n \rightarrow +\infty}|\nabla w_{n}|_{2}^{2}<\frac{2\pi}{\gamma_{0}}-a^{2}-b^{2}$ for all $\mu \geq \mu_{1}$, Corollary \ref{Conver}  together with the equality below
\begin{equation*}
	-\lambda_{1,n}| u_{n}|_{2}^{2}-\lambda_{2,n}| v_{n}|_{2}^{2}=|\nabla w_n|^{2}_{2}-\int_{\mathbb{R}^2}\nabla H(w_{n})\cdot w_{n}\,dx+o_n(1),
\end{equation*}
lead to
\begin{equation}  \label{m2}
	-\lambda_{1,n}a^{2}-\lambda_{2,n}b^{2}=|\nabla w_{n}|_{2}^{2}+o_n(1).
\end{equation}
From this,
\begin{eqnarray*}
    \begin{aligned}\displaystyle
0 \geq \limsup_{n \rightarrow +\infty} (-\lambda_{1,n}a^{2}-\lambda_{2,n}b^{2})= \limsup_{n \rightarrow +\infty} |\nabla w_{n}|^{2}_{2} \geq \liminf_{n \rightarrow +\infty} |\nabla w_{n}|^{2}_{2}\geq 0,
    \end{aligned}
\end{eqnarray*}
then $|\nabla u_n|^{2}_{2} \to 0$, which is absurd, because $m_{\mu}(a,b)>0$.

Now, we obtain that $w_{n}\rightharpoonup w:=(u,v)$, where $u,v\neq 0$. Then, we show the strong convergence that $(u_{n},v_{n})\rightarrow(u,v)$ in $E_{rad}$. The proof is divided into two steps.

{\bf Step 1.} We show that $\lim\limits_{n\rightarrow\infty} |\nabla u_{n}|_{2}^{2}=|\nabla u|_{2}^{2}$ and $\lim\limits_{n\rightarrow\infty}|\nabla v_{n}|_{2}^{2}=|\nabla v|_{2}^{2}$.

By Corollary $\ref{Conver}$, we obtain
\begin{equation}\label{limittb}
    \begin{aligned}\displaystyle
	H(w_n) \to H(w) \quad \mbox{and} \quad \nabla H(w_{n})\cdot w_{n} \to \nabla H(w)\cdot w \,\, \mbox{in} \,\, L^{1}(\mathbb{R}^2).
    \end{aligned}
\end{equation}	
Then, by \eqref{Pohozaev} and $\eqref{limittb}$, together with the weak convergence
\begin{eqnarray*}
    \begin{aligned}\displaystyle
P(u_{n},v_{n})-P(u,v)=o_{n}(1),
    \end{aligned}
\end{eqnarray*}
we deduce
\begin{eqnarray*}
    \begin{aligned}\displaystyle
\lim\limits_{n\rightarrow\infty}|\nabla {u}_{n}|_{2}^{2}-|\nabla {u}|_{2}^{2}+\lim\limits_{n\rightarrow\infty}|\nabla {v}_{n}|_{2}^{2}-|\nabla {v}|_{2}^{2}=0.
    \end{aligned}
\end{eqnarray*}
By Fatou Lemma, we have
\begin{eqnarray*}
    \begin{aligned}\displaystyle
|\nabla {u}|_{2}^{2}\leq\liminf\limits_{n\rightarrow\infty}|\nabla {u}_{n}|_{2}^{2}\leq\lim\limits_{n\rightarrow\infty}|\nabla {u}_{n}|_{2}^{2} \ \mbox{and} \ |\nabla {v}|_{2}^{2}\leq\liminf\limits_{n\rightarrow\infty}|\nabla {v}_{n}|_{2}^{2}\leq\lim\limits_{n\rightarrow\infty}|\nabla {v}_{n}|_{2}^{2},
    \end{aligned}
\end{eqnarray*}
which implies that
\begin{eqnarray*}
    \begin{aligned}\displaystyle
\lim\limits_{n\rightarrow\infty}|\nabla u_{n}|_{2}^{2}=|\nabla u|_{2}^{2} \ \ \mbox{and} \ \ \lim\limits_{n\rightarrow\infty}|\nabla v_{n}|_{2}^{2}=|\nabla v|_{2}^{2}.
    \end{aligned}
\end{eqnarray*}

{\bf Step 2.} We show that $|u|_{2}=a$, $|v|_{2}=b$.

By Lemma \ref{mu}, we have $\lambda_{1}>0$ and $\lambda_{2}>0$. Combining \eqref{mu2}, \eqref{EQ2}, \eqref{PEQ2} and $P(u,v)=0$, we obtain
\begin{eqnarray*}
    \begin{aligned}\displaystyle
    \lambda_{1}a^{2}+\lambda_{2}b^{2}=\lambda_{1}|u|_{2}^{2}+\lambda_{2}|v|_{2}^{2},
    \end{aligned}
\end{eqnarray*}
Then, $0<|u|_{2}\leq a$, and $0<|v|_{2}\leq b$ implying that $|u|_{2}=a$ and $|v|_{2}=b$. Thus, $w_{n}\rightarrow w$ in $E_{rad}$. Finally, by Lemma \ref{fin}, we obtain that $(u,v)$ is a normalized ground state solution of \eqref{aa}.
\qed

\smallskip

\smallskip
\smallskip

\noindent{\bfseries Acknowledgements:}
The research has been supported by National Natural Science Foundation of China 11971392, Natural Science
Foundation of Chongqing, China cstc2019jcyjjqX0022.


\begin{thebibliography}{99}




\bibitem{Alba}
F. S. B. Albuquerque, {\em Nonlinear Schr\"{o}dinger elliptic systems involving exponential critical growth in $\mathbb{R}^{2}$}, Electron. J. Diff. Equ. {\bf 59}(2014), 1-12.

\bibitem{Albb}
F. S. B. Albuquerque, {\em Standing wave solutions for a class of nonhomogeneous systems in dimension two}, Complex Var. Elliptic Equ. {\bf 61}(8)(2016), 1157-1175.


\bibitem{AlvesJi}
C. O. Alves, C. Ji, O. H. Miyagaki, {\em Normalized solutions for a Schr\"{o}dinger equation with critical growth in $\mathbb{R}^{N}$}, Calc. Var. Partial Differ. Equ. {\bf 61}(1)(2022), 1-24.

\bibitem{BartschJeanjean}
T. Bartsch, L. Jeanjean, {\em Normalized solutions for nonlinear Schr\"{o}dinger systems}, Proc. Roy. Soc. Edinburgh Sect. A {\bf 148}(2)(2018), 225-242.

\bibitem{BartschJeanjeanSaove}
T. Bartsch, L. Jeanjean, N. Soave, {\em Normalized solutions for a system of coupled cubic Schr\"odinger equations on $\mathbb{R}^3$}, J. Math. Pures Appl. {\bf 106}(4)(2016), 583-614.

\bibitem{BartschLi}
T. Bartsch, H. W. Li, W. M. Zou, {\em Existence and asymptotic behavior of normalized ground states for Sobolev critical Schr\"{o}dinger systems}, arXiv preprint arXiv:2204.10634, 2022.

\bibitem{BartschSaovea}
T. Bartsch, N. Soave, {\em  A natural constraint approach to normalized solutions of nonlinear Schr\"{o}dinger equations and systems}, J. Funct. Anal. {\bf 272}(12)(2017), 4998-5037.

\bibitem{BartschSaoveb}
T. Bartsch, N. Soave, {\em Multiple normalized solutions for a competing system of Schr\"odinger equations}, Calc. Var. Partial Differ. Equ. {\bf 58}(1)(2019), 1-24.

\bibitem{BartschZhong}
T. Bartsch, X. X. Zhong, W. M. Zou, {\em Normalized solutions for a coupled Schr\"{o}dinger system}, Math. Ann. {\bf 380}(3-4)(2021), 1713-1740.


\bibitem{Cao}
D. M. Cao, {\em Nontrivial solution of semilinear elliptic equation with critical exponent in $\mathbb{R}^2$}, Comm. Partial Differ. Equ. {\bf 17}(3-4)(1992), 407-435.

\bibitem{Cas}
D. Cassani, H. Tavares, J. J. Zhang, {\em Bose fluids and positive solutions to weakly coupled systems with critical growth in dimension two}, J. Differ. Equ. {\bf 269}(3)(2020), 2328-2385.

\bibitem{Chen}
S. T. Chen, A. Fiscella, P. Pucci, X. H. Tang, {\em Coupled elliptic systems in $\mathbb{R}^{N}$ with $(p, N)$ Laplacian and critical exponential nonlinearities}, Nonlinear Anal. {\bf 201}(2020), 1-14.

\bibitem{DOJM}
J. M. do \'{O}, J. C. de Albuquerque, {\em Positive ground state of coupled systems of Schr\"{o}dinger equations in $\mathbb{R}^{2}$ involving critical exponential growth}, Math. Methods Appl. Sci. {\bf 40}(18)(2017), 6864-6879.

\bibitem{DOJMa}
J. M. do \'{O}, J. C. de Albuquerque, {\em On coupled systems of nonlinear Schr\"{o}dinger equations with critical exponential growth}, Appl. Anal. {\bf 97}(6)(2018), 1000-1015.

\bibitem{DOGI}
J. M. do \'{O}, J. Giacomoni, P. K. Mishra, {\em Nonautonomous fractional Hamiltonian system with critical exponential growth}, Nonlinear Differ. Equ. Appl. {\bf26}(4)(2019), 1-25.

\bibitem{Gou}
T. X. Gou, L. Jeanjean, {\em Multiple positive normalized solutions for nonlinear Schr\"odinger systems}, Nonlinearity {\bf 31}(5) (2018), 2319-2345.


\bibitem{jeanjean1}
L. Jeanjean, {\em Existence of solutions with prescribed norm for semilinear elliptic equations}, Nonlinear Anal. {\bf 28}(10)(1997), 1633-1659.


\bibitem{Li}
H. W. Li, W. M. Zou, {\em Normalized ground states for semilinear elliptic systems with critical and subcritical nonlinearities}, J. Fixed Point Theory Appl. {\bf 23}(3)(2021), 1-30.

\bibitem{Moser}
J. Moser, {\em A sharp form of an inequality by N. Trudinger}, Ind. Univ. Math. J. {\bf 20}(1971), 1077-1092.

\bibitem{Thin}
T. V. Nguyen, {\em Existence of solution to singular Schr\"{o}dinger systems involving the fractional p-Laplacian with Trudinger-Moser nonlinearity in $\mathbb{R}^{N}$}, Math. Methods Appl. Sci. {\bf 44}(8)(2021), 6540-6570.

\bibitem{SOAVE1}
N. Soave, {\em Normalized ground states for the NLS equation with combined nonlinearities}, J. Differ. Equ. {\bf 269}(9)(2020), 6941-6987.

\bibitem{SOAVE}
N. Soave, {\em Normalized ground states for the NLS equation with combined nonlinearities: the Sobolev critical case}, J. Funct. Anal. {\bf 279}(6)(2020), 1-43.

\bibitem{Trudinger}
N. S. Trudinger, {\em On imbeddings into Orlicz spaces and some applications}, J. Math. Mech. {\bf 17}(1967), 473-483.

\bibitem{Willem}
M. Willem, {\em  Minimax Theorems}, Birkh\"{a}user 1996.


\end{thebibliography}
\end{document}